\documentclass[11pt]{amsart}

\usepackage{amsfonts}
\usepackage{amsmath}
\newtheorem{theorem}{\bf Theorem}

\def\R{\mathbb{R}} 
\def\K{\mathbf{K}} 
\def\N{\mathbb{N}} 
 
\def\Si{\mathbf{\Sigma^2}}
\begin{document}

\title{The moment problem with bounded density}

\author{Jean B. Lasserre}

\address{LAAS-CNRS and Institute of Mathematics,
LAAS 7 Avenue du Colonel Roche, 31077 Toulouse C\'edex 4 (France)\\
url: www.laas.fr/$\sim$lasserre}
\email{lasserre@laas.fr}



\begin{abstract}
Let $\mu$ be a given Borel measure on $\K\subseteq\R^n$ and 
let $y=(y_\alpha)$, $\alpha\in\N^n$, be a given sequence.
We provide several conditions linking $y$ and
the moment sequence $z=(z_\alpha)$ of $\mu$, for $y$ to be the moment sequence of a Borel measure
$\nu$ on $\K$ which is absolutely continuous with respect to $\mu$ and such that
its density is in  $L_\infty(\K,\mu)$.

The conditions are necessary and sufficient if $\K$ is a compact basic semi-algebraic set, 
and sufficient if $\K\equiv\R^n$. Moreover,  arbitrary finitely many of these conditions can be checked by solving either a semidefinite program or a linear program with a single variable!
\end{abstract}

\maketitle

\section{Introduction}
We consider the following moment problem: Let $\mu$ be a given Borel measure 
$\mu$ on $\K\subseteq\R^n$ and 
let $y=(y_\alpha)$, $\alpha\in\N^n$, be a given sequence. Under what conditions on $y$, 
do we have
\begin{equation}
\label{density}
y_\alpha\,=\,\int_\K\, X^\alpha \,h\,d\mu,\qquad \forall \alpha\in\N^n,
\end{equation}
for some $0\leq h\in L_\infty(\K,\mu)$? 
(The notation $X^\alpha$ stands for the monomial $\prod_iX_i^{\alpha_i}$.)
In Knill \cite{knill}, the measure $d\nu:=hd\mu$ is said to be {\it uniformly absolutely continuous} with respect to $\mu$.

This moment problem (\ref{density}) {\it with density}, 
initially studied by Hausdorff for the interval 
$[0,1]$ of the real line, is a particular case of the 
general $\K$-moment problem where for a given sequence $y$, one only asks for 
existence of some Borel representing measure $\nu$ on $\K$ (not necessarily with a density with respect to some specified measure $\mu$). 
For more details on the general $\K$-moment problem, the interested reader is referred to e.g. Akhiezer \cite{a}, Beg \cite{berg} ,
Curto and Fialkow \cite{curto}, Shohat and Tamarkin \cite{shohat}, Simon \cite{simon}
and many references therein.

The moment problem {\it with density}  is useful in e.g. probability theory, mathematical physics, statistical mechanics, potential theory and dynamical systems, and an interesting discussion on
its potential applications can be found in e.g. Knill \cite{knill}.\\

\noindent
{\bf Contribution.} Our contribution is twofold as we consider the case where $\K\subset\R^n$ is a 
compact basic semi-algebraic set and the case where $\K\equiv\R^n$. 

$\bullet$ In the former case we provide three sets of {\it necessary and sufficient conditions} relating  $y$ with the moment sequence $z$ of $\mu$.
The first two sets of conditions are stated in terms of {\it positive definiteness} of 
so-called moment and localizing matrices, with entries linear in $y$ and $z$, whereas the third 
set of conditions is stated in terms of {\it linear} inequalities linking the moments of
$\mu$ and $\nu$.  The first two  sets of (positive semidefiniteness)
conditions are established thanks to powerful
moment results which are the dual sides of Schm\"udgen's Positivstellensatz \cite{schmudgen},
and its Putinar \cite{putinar} and Jacobi and Prestel \cite{jacobi} refinement.
Checking arbitrary finitely many of these conditions
reduces to solving a {\it semidefinite program} with only one variable! A semidefinite program is
a convex optimization problem for which efficient public softwares are now available.
For more details on the theory and applications of semidefinite programming, 
the interested reader is referred to e.g. Vandenberghe and Boyd \cite{boyd}.

The third set of (linear) conditions is established thanks to
a result by Vasilescu \cite{vasilescu} (see also Krivine \cite{krivine1,krivine2})
which generalizes to basic semi-algebraic sets a previous result of Handelman \cite{handelman}
for polytopes (the latter being itself a generalization to polytopes of the famous Hausdorff 
moment conditions on the unit box $[0,1]^n$). Checking (arbitrary) finitely many of these conditions
reduces to checking whether some interval of the real line is empty.

This third set of linear conditions generalizes to compact basic
semi-algebraic sets a result of Knill \cite{knill}  obtained for the unit box $[0,1]^n$.
This latter result \cite{knill}  states that there is some scalar $M>0$ such
that the classical Haussdorff linear moment conditions
(for a measure to be supported on $[0,1]^n$)
applied to the moment sequence $M\,z$ (for some $M$), {\it dominate} those applied to $y$.

$\bullet$ In the non compact case $\K\equiv\R^n$, we provide a set of {\it sufficient} conditions linking 
$y$ and $z$ and assuming that $z$
satisfies the multivariate Carleman condition. They simply state the existence of a scalar 
$\kappa$ such that
\[0\,\leq\,L_y(f^2)\,\leq\,\kappa\int f^2\,d\mu,\quad\forall f\in\R[X],\]
where $L_y:\R[X]\to\R$ is the linear functional
\[f\,(=\sum_\alpha f_\alpha X^\alpha)\,\mapsto\, L_y(f)\,=\,\sum_\alpha f_\alpha y_\alpha.\]
Notice that these conditions are very simple as they require
that both moment matrices associated with the sequences
$y$ and $\kappa z-y$ are positive semidefinite.
Crucial for establishing the result is the generalized multivariate Carleman condition 
of Nussbaum \cite{nussbaum}.

\section{Notation, definitions and preliminary results}

For a closed set $\K\subset\R^n$, let $\mathcal{B}(\K)$ denote its associated Borel $\sigma$-algebra, and given a Borel measure $\mu$ on $\K$, let $L_\infty(\K,\mu)$ be the standard Banach space of
all bounded measurable functions on $\K$ (except possibly on a subset of $\mu$-measure zero).

Let $\R[X]$ be the ring of real polynomials in the variables $X=(X_1,\ldots,X_n)$, and
let $\Si\subset\R[X]$ be its subset of sums of squares.
With $\Vert .\Vert$ being the euclidean norm of $\R^n$, the notation
$\Vert X\Vert^2$ stands for the quadratic polynomial $\sum_{i=1}^n X_i^2$.
For a real symmetric matrix $A$, the notation $A\succeq0$ satnds for $A$ is positive semidefinite.

Given a family $\{g_j\}_{j=1}^m\subset\R[X]$,
denote by:

$\bullet$  $g_J\in\R[X]$, the polynomial
$g_J(X):=\prod_{j\in J}g_j(X)$, for every $J\subseteq\{1,\ldots,m\}$, 
and  with the convention $g_\emptyset\equiv 1$.

$\bullet$  $Q(g)\subset\R[X]$, the {\it quadratic module} generated by the $g_j$'s, that is,
\[Q(g)\,:=\,\{\:u_0+\sum_{j=1}^mu_j\,g_j\quad \vert\quad u_j\in\Si\quad \forall j=0,1,
\ldots,m\}.\]
Given a sequence $y=(y_\alpha)$ indexed in the canonical basis $(X^\alpha)$ of $\R[X]$, let $L_y:\R[X]\to \R$ be the linear functional
\[f\,(=\sum_{\alpha\in\N^n}\, f_\alpha X^\alpha)\,\mapsto\,L_y(f)\,:=\,\sum_{\alpha\in \N^n}\,f_\alpha y_\alpha.\]
Given a Borel set $\K\subseteq\R^n$, a sequence $y=(y_\alpha)$ has a {\it representing measure}
 $\mu$ on $\K\subset\R^n$ if
 \[y_\alpha\,=\,\int_\K X^\alpha\,d\mu,\qquad \forall\,\alpha\in\N^n.\]
A Borel measure $\mu$ on $\R^n$ 
with all its moments finite is said to be {\it moment determinate} if
whenever
\[\int X^\alpha \,d\mu \,=\,\int X^\alpha \,d\nu,\qquad \forall \,\alpha\in\N^n,\]
holds for some measure $\nu$, then $\nu=\mu$.

When $\K\equiv\R^n$ and if
\begin{equation}
\label{carleman}
\sum_{k=1}^\infty L_y(X_i^{2k})^{-1/2k}\,=\,\infty,
\end{equation}
then  by a result from Nussbaum \cite{nussbaum},  $y$ has a unique representing measure $\mu$; see e.g. Berg \cite{berg}. Condition (\ref{carleman}) is called the generalized Carleman condition as 
an extension to the multivariate case of Carleman's condition for the univariate case.

We next recall several important results concerning necessary and sufficient conditions for
a sequence $y=(y_\alpha)$ to have a representing measure on  a compact basic semi-algebraic set $\K\subset\R^n$.

So, let $\K\subset\R^n$ be a basic semi-algebraic set defined by:
\begin{equation}
\label{setk}
\K:=\{x\in\R^n\:\vert \quad g_j(x)\geq0,\quad\forall j=1,\ldots,m\},
\end{equation}
for some polynomials $\{g_j\}\subset\R[X]$.

\begin{theorem}
\label{rappel}
Let $\K\subset\R^n$ be compact and defined as in (\ref{setk}). 

{\rm (a)} A given sequence $y=(y_\alpha)$ has a representing
measure $\mu$ on $\K$ if and only if:
\begin{equation}
\label{rappel-1}
L_y(f^2\,g_J)\,\geq\,0\qquad\forall \,f\in\R[X];\quad \forall J\subseteq\{1,\ldots,m\}.
\end{equation}

{\rm (b)} Assume that  for some $M>0$, the quadratic polynomial $M-\Vert X\Vert^2$ is in $Q(g)$. 
A given sequence $y=(y_\alpha)$ has a representing measure $\mu$ on $\K$ if and only if:
\begin{equation}
\label{rappel-2}
L_y(f^2\,g_j)\,\geq\,0\qquad\forall \,f\in\R[X];\quad \forall j=0,1,\ldots,m
\end{equation}
(with the convention $g_0\equiv 1$).
\end{theorem}
Theorem \ref{rappel}(a) is the dual part of Schm\"udgen's Positivstellensatz \cite{schmudgen} whereas
its refinement (b) is due to Putinar \cite{putinar} and Jacobi and Prestel \cite{jacobi}.

If one restricts $f$ to have degree less than, say $r$, then
testing conditions (\ref{rappel-1}) and (\ref{rappel-2}) reduces to 
checking whether some real-valued symmetric matrices are positive semidefinite.
Indeed, given $\theta\in\R[X]$, let $M_r(\theta\,y)$ be the real symmetric matrix with rows and columns indexed 
in the canonical basis $(X^\alpha)$ of $\R[X]$, and defined by 
\[M_r(\theta\,y)(\alpha,\beta)\,=\,L_y(X^{\alpha+\beta}\theta(X)),\qquad
\forall \alpha,\beta\in\N^n\:\mbox{ with }\vert\alpha\vert,\vert\beta\vert\leq r\]
(where for every $\alpha\in\N^n$,  $\vert\alpha\vert:=\sum_i\alpha_i$).
$M_r(\theta y)$ is called the {\it localizing} matrix associated with $y$ and $\theta$, and when $\theta\equiv 1$, $M_r(y)$ is called the {\it moment} matrix associated with  $y$.

Then testing condition (\ref{rappel-1}) (resp. (\ref{rappel-2})) for all $f\in\R[X]$ 
with ${\rm deg}\,f\leq r$, is equivalent to checking whether
$M_r(g_J\,y)\succeq0$ for all $J\subseteq\{1,\ldots,m\}$ 
(resp. $M_r(g_j\,y)\succeq0$ for all $j=0,\ldots,m$), i.e., computing 
the smallest eigenvalue of $2^m$ (resp. $m+1$) real 
symmetric matrices. If $y$ is taken as unknown, then 
the LMIs (Linear Matrix Inequalities) $M_r(g_J\,y)\succeq0$ (or 
$M_r(g_j\,y)\succeq0$) on $y$ define the set of sequences $y$ 
which have a representing measure $\mu$ on $\K$. Importantly, this set is convex.

Finally, given $\alpha\in\N^n$, let $g^\alpha,(1-g)^\alpha\in\R[X]$ be the polynomials
\[g^\alpha(X)\,:=\,\prod_{j=1}^m\,g_j(X)^{\alpha_j};\quad
(1-g)^\alpha(X)\,:=\,\prod_{j=1}^m\,(1-g_j(X))^{\alpha_j}.\]
\begin{theorem}
\label{rappel3}
Let $\K\subset\R^n$ be as in (\ref{setk}), and assume that the $g_j$'s are normalized, i.e., 
$0\leq g_j\leq1$ on $\K$ for all
$j=1,\ldots,m$, and that the family $(0,1,\{g_j\})$ generates the algebra $\R[X]$; equivalently
$\R[X]=\R[g_1,\ldots,g_m]$.

A given sequence $y=(y_\alpha)$ has a representing
measure $\mu$ on $\K$ if and only if:
\begin{equation}
\label{rappel3-1}
L_y(g^\alpha(1-g)^\beta)\,\geq\,0\qquad\forall \,\alpha,\beta\in\N^n.
\end{equation}
\end{theorem}
Theorem \ref{rappel3} is from Vasilescu \cite{vasilescu} (also
hidden in a previous result of Krivine \cite{krivine1,krivine2}).
Notice that testing whether (\ref{rappel3-1}) holds for all $\alpha,\beta\in\N^n$ with
$\vert\alpha\vert,\vert\beta\vert\leq r$, reduces to checking
whether finitely many inequalities are satisfied. If again $y$ is taken as unknown, then
the {\it linear} inequalities (\ref{rappel3-1}) on $y$ are another representation of the convex
set of feasible sequences $y$ which have a representing measure on $\K$.

So this set can be represented either via LMIs through conditions (\ref{rappel-1}) (or (\ref{rappel-2})),
or via linear inequalities through conditions (\ref{rappel3-1}).

\section{Main result}

Let $\mu$ be a finite Borel measure on a Borel set $\K\subseteq\R^n$. We want to find a set of necessary and/or 
sufficient conditions on a sequence $y=(y_\alpha)$ to have a representing measure $\nu$ on $\K$ that is {\it uniformly 
absolutely continuous} with respect to $\mu$, that is, $\nu$ is absolutely continuous
w.r.t. $\mu$ (denoted $\nu\ll\mu$) and its Radon-Nikodym
derivative $h$ belongs to $L_\infty(\K,\mu)$. We first concentrate 
on the case where $\K\subset\R^n$ is basic compact semi-algebraic set as defined in (\ref{setk}), and then consider the non compact case $\K\equiv\R^n$.
\subsection{The compact case}
\begin{theorem}
\label{thmain}
Let $\K\subset\R^n$ be compact and defined as in (\ref{setk}).
Let $z=(z_\alpha)$ be the moment 
sequence of a finite Borel measures $\mu$ on $\K$. 

{\rm (a)} A sequence $y=(y_\alpha)$ has a representing Borel measure 
$\nu$ on $\K$, uniformly absolutely continuous with respect to $\mu$, if and only if there is some scalar
$\kappa$ such that for every
$J\subseteq\{1,\ldots,m\}$:
\begin{equation}
\label{thmain-1}
0\,\leq L_y(f^2g_J)\,\leq\,\kappa\,L_z(f^2g_J)\qquad\forall \,f\in\R[X].
\end{equation}

{\rm (b)} In addition, if the polynomial $N-\Vert X\Vert^2$ belongs to the quadratic module $Q(g)$ then
one may replace (\ref{thmain-1}) with the weaker condition
\begin{equation}
\label{thmain-2}
0\,\leq L_y(f^2g_j)\,\leq\,\kappa\,L_z(f^2g_j)\quad\forall \,f\in\R[X];\quad j=0,\ldots,m
\end{equation}
(with the convention $g_0\equiv 1$).

{\rm (c)} Suppose that the $g_j$'s are normalized so that
\[0\leq g_j\leq 1\quad\mbox{on }\quad\K,\quad \forall\,j=1,\ldots,m,\]
and that the family $(0,1,\{g_j\})$ generates the algebra $\R[X]$.

A sequence $y=(y_\alpha)$ has a representing Borel measure 
$\nu$ on $\K$, uniformly absolutely continuous with respect to $\mu$, if and only if there is some scalar
$\kappa$ such that for every $\alpha,\beta\in\N^m$:
\begin{equation}
\label{thmain-3}
0\,\leq L_y(g^\alpha (1-g)^\beta)\,\leq\,\kappa\,L_z(g^\alpha (1-g)^\beta).
\end{equation}
\end{theorem}
\begin{proof}
(a). The {\it only if} part:
Let  $d\nu =hd\mu$ for some $0\leq h\in L_\infty(\K,\mu)$, and let
$\kappa:=\Vert h\Vert_\infty$. Observe that
$g_J\geq0$ on $\K$ for all $J\subseteq\{1,\ldots,m\}$.
Therefore, for every $J\subseteq\{1,\ldots,m\}$:
\[\int f^2g_J\,d\nu=\int f^2g_Jh\,d\mu\,\leq\,\kappa\int f^2g_J\,d\mu,\qquad\forall \,f\in\R[X],\]
and so (\ref{thmain-1}) is satisfied.

The {\it if part}: Let $y$ and $z$ be such that (\ref{thmain-1}) holds true.
Then by Theorem \ref{rappel}(a),
$y$ has a representing finite Borel measure $\nu$ on $\K$.
In addition, let $\gamma=(\gamma_\alpha)$ with $\gamma_\alpha
:=\kappa z_\alpha-y_\alpha$ for all $\alpha\in\N^n$.
From (\ref{thmain-1}), one has
\[L_\gamma(f^2g_J)\geq0,\qquad\forall \,f\in\R[X];\quad \forall J\subseteq\{1,\ldots,m\}.\]
By Theorem \ref{rappel}(a) again, $\gamma$ has a representing
finite Borel measure $\psi$ on $\K$. Moreover, from the definition of the sequence 
$\gamma$
\[\int f \,d(\psi+\nu)\,=\,\int f \,\kappa\,d\mu,\qquad\forall\,f\in\R[X],\]
and therefore, as measures on compact sets are moment determinate,
$\psi+\nu=\kappa\mu$. Hence $\kappa\mu\geq\nu$
which shows that $\nu\ll\mu$. Finally, write $d\nu =h\,d\mu$ for some $0\leq h\in L_1(\K,\mu)$.
From $\nu\leq\kappa\mu$ one obtains
\[\int_A (h-\kappa)\,d\mu \qquad\forall\,A\,\in\,\mathcal{B}(\K),\]
and so $0\leq h\leq \kappa$, $\mu$-almost everywhere on $\K$. Equivalently,
$\Vert h\Vert_\infty\leq\kappa$, the desired result.

(b) The proof is a verbatim copy of (a), except that instead of invoking
Theorem \ref{rappel}(a), we now invoke Theorem \ref{rappel}(b), its Putinar's refinement.

(c) The {\it only if} part:
Let  $d\nu =hd\mu$ for some $h\in L_\infty(\K,\mu)$, and let
$\kappa:=\Vert h\Vert_\infty$. As $g\geq0$ and $(1-g)\geq0$ on $\K$, 
\[\int g^\alpha (1-g)^\beta\,d\nu=\int g^\alpha (1-g)^\beta h\,d\mu\,\leq\,\kappa
\int g^\alpha (1-g)^\beta \,d\mu,\]
for all $\alpha,\beta\in\N^m$, and so (\ref{thmain-3}) is satisfied.

The {\it if part}: Let $y$ and $z$ be such that (\ref{thmain-3}) holds true.
By Theorem \ref{rappel3}, $y$ has a representing
Borel measure $\nu$ on $\K$.
In addition, let $\gamma=(\gamma_\alpha)$ with $\gamma_\alpha
:=\kappa z_\alpha-y_\alpha$ for all $\alpha\in\N^n$.
From (\ref{thmain-3}), the condition (\ref{rappel3-1}) holds true
with $\gamma$ in lieu of $y$, and therefore, by Theorem \ref{rappel3},
$\gamma$ has a representing Borel measure $\psi$ on $\K$. Moreover, from the definition of 
$\gamma$, one has
\[\int f \,d(\psi+\nu)\,=\,\int f \,\kappa\,d\mu,\qquad\forall\,f\in\R[X],\]
and the rest of the proof is the same as that of (a).
\end{proof}
Recalling the definition of localizing matrices,
notice that checking (\ref{thmain-1}) for all $f\in\R[X]$ with ${\rm deg}\,f\leq r$, reduces to checking whether
\begin{equation}
\label{1}
0\,\preceq\,M_r(g_J\,y)\quad\mbox{and}\quad
M_r(g_J\,y)\preceq\,\kappa\,M_r(g_J\,z)\qquad\forall J\subseteq
\{1,\ldots,m\},
\end{equation}
whereas checking (\ref{thmain-2}) for all $f\in\R[X]$ with ${\rm deg}\,f\leq r$ reduces to checking whether
\begin{equation}
\label{2}
0\,\preceq\,M_r(g_j\,y)\quad\mbox{and}\quad
M_r(g_j\,y)\preceq\,\kappa\,M_r(g_j\,z)\qquad\forall j=0,\ldots,m.
\end{equation}
While checking the first LMIs in the left of both
both (\ref{1}) and (\ref{2}) is just evaluating the smallest eigenvalue of
the $2^m$ symmetric matrices $M_r(g_J\,y)$, $J\subseteq\{1,\ldots,m\}$ (resp.
the $m+1$ symmetric matrices $M_r(g_j\,y)$, $j=0,1,\ldots,m$), checking
the second LMIs in both (\ref{1}) and (\ref{2})
reduces to solving  a semidefinite program with only {\it one} variable $\kappa$,
and $2^{m}$ (resp. $m+1$)  LMIs, a convex optimization problem for which efficient
public softwares are now available. \\

On the other hand, checking (\ref{thmain-3}) or all $f\in\R[X]$ with ${\rm deg}\,f\leq r$, reduces 
to solving a linear program with only one variable also, i.e., checking whether some interval is nonempty.

Theorem \ref{thmain}(c) is the analogue for compact basic semi-agebraic sets of
Corollary 3.1 in Knill \cite{knill} for the unit box $[0,1]^n$.

We next obtain a sufficient condition for the non-compact case.

\subsection{The non-compact case}
Let $\mu$ be a finite Borel measure on $\R^n$. We want to find a 
sufficient condition on a sequence $y=(y_\alpha)$ to be the moment sequence
of a measure $\nu$ on $\R^n$ that is {\it uniformly 
absolutely continuous} with respect to $\mu$.

\begin{theorem}
\label{th2}
Let $z=(z_\alpha)$ be the moment 
sequence of a finite Borel measure $\mu$ on $\R^n$, and assume that
$z$ satisfies the generalized Carleman condition, i.e.,
\begin{equation}
\label{th2-1}
\sum_{k=1}^\infty L_z(X_i^{2k})^{-1/2k}\,=\,\infty.
\end{equation}
Then a sequence $y=(y_\alpha)$ has a representing Borel measure 
$\nu$ on $\R^n$, uniformly absolutely continuous with respect to $\mu$,  if 
there exists a scalar $0<\kappa$ such that for all $i=1,\ldots,n$:
\begin{equation}
\label{th2-3}
0\,\leq\,L_y(f^2)\,\leq\,\kappa\,L_z(f^2)\qquad\forall \,f\in\R[X].
\end{equation}
\end{theorem}
\begin{proof}
For every $i=1,\ldots,n$,  (\ref{th2-3}) with $f^2=X_i^{2k}$ yields
\[L_y(X_i^{2k})^{-1/2k}\,\geq\,  \kappa^{-1/2k}\,L_z(X_i^{2k})^{-1/2k},
\qquad\forall\,k=0,1,\ldots,\]
and so, using (\ref{th2-1}), one obtains
\[\sum_{k=1}^\infty \,L_y(X_i^{2k})^{-1/2k}\,\geq\,
\sum_{k=1}^\infty \,\kappa^{-1/2k}\,L_z(X_i^{2k})^{-1/2k}\,=\,+\infty,\]
for every $i=1,\ldots,n$, i.e., the generalized Carleman 
condition (\ref{carleman}) holds for the sequence $y$. Combining this with the first inequality
in (\ref{th2-3})  yields that 
$y$ has a unique representing
Borel measure $\nu$ on $\R^n$.
It remains to prove that $\nu\ll\mu$ and its density $h$ is
in $L_\infty(\R^n,\mu)$. 

Let $\gamma=(\gamma_\alpha)$ with $\gamma_\alpha
:=\kappa z_\alpha-y_\alpha$ for all $\alpha\in\N^n$. 
Then the second inequality in (\ref{th2-3}) yields
\begin{equation}
\label{testt}
L_\gamma(f^2)\geq0\qquad\forall f\in\R[X].\end{equation}
Next, observe that from (\ref{th2-3}), for every $i=1,\ldots,n$, and every $k=0,1,\ldots$,
\[L_\gamma(X_i^{2k})\,\leq\,  \kappa\,L_z(X_i^{2k}),\]
which implies
\begin{equation}
\label{test}
L_\gamma(X_i^{2k})^{-1/2k}\,\geq\,  \kappa^{-1/2k}\,L_z(X_i^{2k})^{-1/2k},\end{equation}
and so, for every $i=1,\ldots,n$, 
\[\sum_{k=1}^\infty \,L_\gamma(X_i^{2k})^{-1/2k}\,\geq\,
\sum_{k=1}^\infty \,\kappa^{-1/2k}\,L_z(X_i^{2k})^{-1/2k}\,=\,+\infty,\]
i.e., $\gamma$ satisfies the generalized Carleman condition. In view of (\ref{testt}), $\gamma$ has a (unique) representing measure $\psi$ on $\R^n$. Next, 
from the definition of $\gamma$, one has
\[\int f\,d(\psi+\nu)\,=\,\kappa\,\int f\,d\mu,\qquad\forall\,f\in\R[X].\]
But as $\mu$ (and so $\kappa\mu$) satisfies Carleman condition
(\ref{th2-1}), $\kappa\mu$ is moment determinate and therefore,  $\kappa\mu=\nu+\psi$.

Hence $\nu\ll\mu$ follows from $\nu\leq\kappa\mu$. Finally, writing
$d\nu=hd\mu$ for some nonnegative $h\in L_1(\R^n,\mu)$, and using $\nu\leq\kappa\mu$, one obtains
\[\int_A (h-\kappa)\,d\mu \,\leq\,0\qquad\forall\,A\,\in\,\mathcal{B}(\R^n),\]
and so $0\leq h\leq \kappa$, $\mu$-almost everywhere on $\R^n$. Equivalently,
$\Vert h\Vert_\infty\leq\kappa$, the desired result.
\end{proof}
Observe that condition (\ref{th2-3}) is extremely simple as it 
is equivalent to stating that
\[\kappa\,M_r(z)\,\succeq\,M_r(y)\,\succeq\,0,\qquad\forall r=0,1,\ldots\]
where $M_r(z)$ (resp. $M_r(y)$) is the $r$-moment matrix associated 
with the sequence $z$ (resp. $y$).
Finally, a sufficient condition for (\ref{th2-1}) to hold is 
\[\int \exp{\vert x_i\vert}\,d\mu\,<\,\infty,\qquad \forall i=1,\ldots,n.\]

\section*{Acknowledgements}
The author acknowledges financial support by the (french) National Research Agency (ANR), under grant NT05-3-41612.

\end{document}